\documentclass{amsart}
\usepackage{amsfonts, amsmath, latexsym, epsfig}
\usepackage[norelsize]{algorithm2e}
\usepackage{amssymb}
\usepackage{color}
\usepackage{epsf}
\usepackage{url}

\newcommand{\RR}{\ensuremath{\mathbb{R}}}

\newcommand{\QQ}{\ensuremath{\mathbb{Q}}}

\newcommand{\ZZ}{\ensuremath{\mathbb{Z}}}

\newtheorem{proposition}{Proposition}
\newtheorem{theorem}{Theorem}
\newtheorem{corollary}{Corollary}
\newtheorem{lemma}{Lemma}

\newtheorem{conjecture}{Conjecture}

\newtheorem{remark}{Remark}
\newtheorem{definition}{Definition}
\def\QuotS#1#2{\leavevmode\kern-.0em\raise.2ex\hbox{$#1$}\kern-.1em/\kern-.1em\lower.25ex\hbox{$#2$}}

\newcommand{\MC}{{\mathcal C}}
\newcommand{\MF}{{\mathcal F}}

\DeclareMathOperator{\Aut}{Aut}

\DeclareMathOperator{\CUTP}{CUTP}
\DeclareMathOperator{\CCUTP}{CCUTP}
\DeclareMathOperator{\GL}{GL}
\DeclareMathOperator{\AGL}{AGL}

\DeclareMathOperator{\conv}{conv}

\DeclareMathOperator{\Stab}{Stab}
\DeclareMathOperator{\Ker}{Ker}

\DeclareMathOperator{\Hyp}{Hyp}
\DeclareMathOperator{\rankperf}{rankperf}

\DeclareMathOperator{\Cst}{Cst}
\DeclareMathOperator{\Lin}{Lin}
\DeclareMathOperator{\Quad}{Quad}
\DeclareMathOperator{\degrk}{degrk}
\DeclareMathOperator{\Vol}{vol}
\DeclareMathOperator{\Aff}{Aff}

\begin{document}

\author{Mathieu Dutour Sikiri\'c}
\address{Mathieu Dutour Sikiri\'c, Rudjer Boskovi\'c Institute, Bijenicka 54, 10000 Zagreb, Croatia, Fax: +385-1-468-0245}
\email{mathieu.dutour@gmail.com}

\thanks{The author gratefully acknowledge support from the Alexander von Humboldt foundation}

\title{The seven dimensional perfect Delaunay polytopes and Delaunay simplices}

\date{}

\maketitle

\begin{abstract}
For a lattice $L$ of $\RR^n$, a sphere $S(c,r)$ of center $c$ and radius $r$
is called {\em empty} if for any $v\in L$ we have $\Vert v - c\Vert \geq r$.
Then the set $S(c,r)\cap L$ is the vertex set of a {\em Delaunay polytope}
$P=\conv(S(c,r)\cap L)$.
A Delaunay polytope is called {\em perfect} if any affine transformation
$\phi$ such that $\phi(P)$ is a Delaunay polytope is necessarily an isometry
of the space composed with an homothety.

Perfect Delaunay polytopes are remarkable structure that exist only
if $n=1$ or $n\geq 6$ and they have shown up recently in covering maxima studies.
Here we give a general algorithm for their enumeration that relies on
the Erdahl cone.
We apply this algorithm in dimension $7$ which allow us to find
that there are only two perfect Delaunay polytopes: $3_{21}$ which
is a Delaunay polytope in the root lattice $\mathsf{E}_7$ and the
Erdahl Rybnikov polytope.

We then use this classification in order to get the list of all types
Delaunay simplices in dimension $7$ and found $11$ types.
\end{abstract}

\section{Introduction}
A lattice $L$ is a set of the form $L=\ZZ v_1 + \dots + \ZZ v_n\subset \RR^n$
with $(v_1, \dots, v_n)$ being independent.
For such $L$ a sphere $S(c,r)$ of center $c$ and radius $r$ is called {\em empty} if for any $v\in L$ we have $\Vert v - c\Vert \geq r$.
A polytope $P$ is called a {\em Delaunay polytope} if it is full-dimensional and if the vertex-set of $P$ is $S(c,r)\cap L$ with $S(c,r)$ an empty sphere.
A Delaunay polytope is called {\em perfect} if any affine transformation
$\phi$ such that $\phi(P)$ is a Delaunay polytope is necessarily an isometry
of the space composed with an homothety.

In \cite{Erdahl1992,DGL92} it was proved that for dimension $n\leq 5$ the only
possible perfect Delaunay polytope is the interval $[0,1]$. Also in \cite{DGL92}
it is proved that the Gosset's polytopes $2_{21}$ and $3_{21}$ which are Delaunay
polytopes of $\mathsf{E}_6$ and $\mathsf{E}_7$ are perfect.
From the construction of infinite sequences in \cite{Dutour2005,Grishukhin2006,InfiniteER2002,SecInfinite}
we know that for any dimension $n\geq 6$ there exist perfect Delaunay polytopes.
In \cite{Dutour2005} for any $n\geq 6$ we define a Delaunay polytope $ED_n$ of a lattice
$LD_n$. The lattice $LD_n$ is formed by lamination over the root lattice
$\mathsf{D}_{n-1}$ and we prove in \cite{InhomogeneousExtremeForms} that $ED_n$
is the unique Delaunay polytope of maximum circumradius of $LD_n$ and compute its
covering density.

In \cite{Hyp7} we proved that the $ED_6 = 2_{21}$ is the unique perfect Delaunay
polytopes in dimension $6$.
This work uses a new approach in order to prove the following theorem:
\begin{theorem}\label{Theorem3_21_ER7}
The $7$-dimensional perfect Delaunay polytopes are the Gosset polytope $ED_7 = 3_{21}$ and
the Erdahl and Rybnikov polytope $ER_7$ \cite{35tope,InfiniteER2002}.
\end{theorem}
\proof See Section \ref{DimSeven_perfectDelaunay}. \qed

Perfect Delaunay polytope are of importance for the theory of Covering Maxima.
A {\em covering maximum} is a lattice $L$ such that its covering density
is reduced if it is perturbed.
In \cite{InhomogeneousExtremeForms} it is proved that a lattice $L$ is a covering
maximum if and only if the Delaunay polytopes of maximum circumradius are perfect
and {\em eutactic} (see \cite{InhomogeneousExtremeForms} for the definition).
This characterization echoes Voronoi's theorem \cite{VoronoiI} for the characterization
of lattices of maximum density in terms of perfection and eutacticity.
In \cite{InhomogeneousExtremeForms} we proved that $LD_n$ is one such covering maxima.
Based on Theorem \ref{Theorem3_21_ER7} and partial enumerations 
in dimension $8$, $9$ and $10$ we state the following conjecture:

\begin{conjecture}\label{UpperBoundConj}
  For each $n\geq 6$, the lattice $LD_n$ defined in \cite{Dutour2005} has
  maximal covering density among all covering maxima.
\end{conjecture}

The Minkowski conjecture \cite[p. 18]{Koksma} on the product of inhomogeneous forms
has inspired a lot of research. Recently, it has been proved for $n\leq 8$ in
\cite{HansGill1,HansGill2,HansGill3,HansGill4} by computational methods based on
Korkine-Zolotarev reduction theory. Other theoretical approaches have been attempted
in \cite{McMullenCurtis,StableLattices} by Dynamical System Theory. In particular the
following theorem is proved in \cite[Corollary 1.3]{StableLattices}:

\begin{theorem}
If Conjecture \ref{UpperBoundConj} holds for a dimension $n\geq 1$ then Minkowski's conjecture
holds for dimension $n$.
\end{theorem}
As a consequence of the work of this paper we have that Minkowski's conjecture is
correct in dimension $7$, thereby confirming \cite{HansGill4}.

We prove Theorem \ref{Theorem3_21_ER7} by using 
the {\em Erdahl cone} which is defined as the set of polynomial functions $f$ of degree
at most $2$ such that $f(x)\geq 0$ for $x\in \ZZ^n$. We already used this cone in
\cite{InhomogeneousExtremeForms} for the study of covering maxima.
We have then to do a kind of dual description computation with the problem that
the number of defining inequalities is infinite, we have no local polyhedrality result
as in the perfect form case (see \cite{bookschurmann} for details) and we are interested
in only a subset of the extreme rays
(see Theorem \ref{ClassifKindExtRay} for the list of possible kinds of extreme rays of the Erdahl cone).

In \cite{Hyp7} we used a different approach, i.e. hypermetrics that
allowed us to find all the $6$-dimensional perfect Delaunay polytopes. But this
approach relied on previous work \cite{Ba,BaRy} on $6$-dimensional Delaunay simplices that
we could not extend easily to dimension $7$. Thus, it appears that the only way
to classify the perfect $7$-dimensional Delaunay is to use the Erdahl cone.
Moreover, we are able to use this classification in order to get the classification
of Delaunay simplices:

\begin{theorem}\label{EnumerationSimplices}
Up to arithmetic equivalence there are $11$ types of $7$-dimensional Delaunay simplices.
The full list is given in Table \ref{TableFundamentalSimplices}.
\end{theorem}
\proof See Section \ref{DimSeven_simplices}. \qed

In contrast to perfect Delaunay polytopes, the lattices simplices of this list (except
the trivial simplex) had not been discovered before. A similar study has been undertaken
in \cite{MinkowskianSublattices} for the set of shortest vectors of lattices. In view of
this work it seems reasonable to think that the classification of Delaunay simplices is
possible in dimension $8$. Equally importantly the classification of perfect Delaunay
polytopes in dimension $8$ could be done and a conjectural list of the $27$ known possibilities
is available in \cite{PerfectDelaunayLowDim}.

\begin{table}
\begin{center}
\caption{Representative of Delaunay simplices in dimension $7$. $e_1, \dots, e_7$ is the standard basis of $\ZZ^7$ and $e_0=0$. $\Vol(S)$ is $n!$ times the Euclidean volume of $S$. $\left\vert\Stab(S)\right\vert$ is the size of the lattice automorphism group preserving $S$. ``Nb interval'' is the number of Delaunay polyhedra of the type $\{0,1\}\times \ZZ^6$ in which $S_i$ is contained}
\label{TableFundamentalSimplices}
\begin{tabular}{||c|c|c|c|c||}
\hline
\hline
$i$ & Representative $S_i$ & $\Vol(S_i)$ & $\left\vert \Stab(S_i)\right\vert$ & Nb interval \\
\hline
1 & $e_0, \dots, e_6,e_7$             & 1 & 40320 & 127\\
2 & $e_0, \dots, e_6,(1,1,1,1,1,1,2)$ & 2 & 40320 & 63\\
3 & $e_0, \dots, e_6,(0,0,1,1,1,1,2)$ & 2 & 1440 & 63\\
4 & $e_0, \dots, e_6,(1,2,2,2,2,2,3)$ & 3 & 5040 & 42\\
5 & $e_0, \dots, e_6,(1,1,1,2,2,2,3)$ & 3 & 1152 & 42\\
6 & $e_0, \dots, e_6,(0,1,1,2,2,2,3)$ & 3 & 240 & 41\\
7 & $e_0, \dots, e_6,(1,1,1,1,1,3,4)$ & 4 & 1440 & 27\\
8 & $e_0, \dots, e_6,(1,1,1,1,2,2,4)$ & 4 & 240 & 31\\
9 & $e_0, \dots, e_6,(1,1,1,2,2,3,4)$ & 4 & 144 & 31\\
10 & $e_0, \dots, e_6,(1,1,3,3,3,4,5)$ & 5 & 72 & 24\\
11 & $e_0, \dots, e_6,(1,1,1,1,2,3,5)$ & 5 & 48 & 24\\
\hline
\hline
\end{tabular}
\end{center}
\end{table}

In Section \ref{Sec:DelaunayPolyhedra} Delaunay polyhedra are considered,
their basic structure and relation to the Erdahl cone are introduced here.
The facial structure of the Erdahl cone is reviewed in Section 
\ref{FacialErdahlCone}, in particular not all extreme rays of the Erdahl
cone are related to Delaunay polyhedra \cite{Erdahl1992}.
We also explain how the hypercube $[0,1]^n$ correspond to the cut polytope
in the Erdahl cone.
Section \ref{GeometryErdahlCone} is not used in later sections.
In it we construct a retraction of the Erdahl cone on the faces defined by
Delaunay polyhedra.
In Section \ref{SectionLtype} we establish the link between the Erdahl cone and
the classic $L$-type theory.
In Section \ref{SectionHypermetric} we do the same for the hypermetric cone.
In Section \ref{SectionConn} we give the connectivity and finiteness results
on which our enumeration algorithm relies.
Then we present in Section \ref{SectionAlgo} our enumeration method,
which is modelled on the Voronoi algorithm for perfect forms
\cite{martinet} and on the adjacency decomposition method  \cite{symsurvey}.
In Section \ref{DimSeven_perfectDelaunay} we give the obtained results in
the classification of $7$-dimensional perfect Delaunay polyhedra.
In Section \ref{DimSeven_simplices} we use this classification
to classify the $7$-dimensional types of Delaunay simplices.

\section{Delaunay polyhedra}\label{Sec:DelaunayPolyhedra}

Denote by $E_2(n)$ the vector space of polynomials of degree at most $2$
on $\RR^n$ and by $\AGL_n(\ZZ)$ the group of affine integral transformations on $\ZZ^n$.
The {\em Erdahl cone} is defined as
\begin{equation*}
Erdahl(n)=\{f\in E_2(n)\mbox{~such~that~}f(x)\geq 0\mbox{~for~}x\in \ZZ^n\}.
\end{equation*}
It is a convex cone of dimension $(n+1)(n+2)/2$ on which the group $\AGL_n(\ZZ)$ acts.
All defining inequalities $f(x)\geq 0$ are equivalent under $\AGL_n(\ZZ)$ and therefore
$Erdahl(n)$ is not polyhedral.

We denote $\cdot$ the standard scalar product on $\RR^n$ defined by $x\cdot y = x^{T} y$.
For a symmetric matrix $A$ and $x\in \RR^n$ we define $A[x] = x^T A x$.
We write any $f\in E_2(n)$ in the form
\begin{equation*}
f(x)=\Cst(f)  + 2\Lin(f) \cdot x + \Quad(f)[x]
\end{equation*}
with $\Cst(f)\in \RR$, $\Lin(f)\in \RR^n$ and $\Quad(f)$ a $n\times n$ symmetric matrix.
We define $S^n$ to be the set of symmetric matrices, $S^n_{>0}$ the set
of positive
definite matrices, $S^{n}_{\geq 0}$ the set of positive semidefinite matrices.
We also define $Erdahl_{>0}(n)$ to be the set of $f\in Erdahl(n)$ with
$\Quad(f)\in S^n_{>0}$.
If $f\in Erdahl(n)$ then $\Quad(f)\in S^n_{\geq 0}$.

A {\em sublattice} of $\ZZ^n$ is a subgroup of $\ZZ^n$. An {\em affine sublattice}
is one of the form $x_0 + L$ with $x_0\in \ZZ^n$ and $L$ a sublattice of $\ZZ^n$.
A lattice $L\subset \ZZ^n$ is called {\em saturated} if $(L\otimes \RR)\cap \ZZ^n = L$.
If $L_1$ and $L_2$ are two sublattices of $\ZZ^n$ we write
$L_1\oplus_{\ZZ} L_2=\ZZ^n$ if $L_1\cap L_2=\{0\}$ and $L_1\cup L_2$
generates $\ZZ^n$ over $\ZZ$. In that case both $L_1$ and $L_2$ are saturated.

\begin{definition}
Let us fix $n\geq 1$ and define:
  
(i) A {\em Delaunay polyhedron} $D$ to be a set of the form $D=P_{L'}(D) + L(D)\subset \ZZ^n$ with 
\begin{itemize}
\item $\conv(P_{L'}(D))$ a Delaunay polytope of an affine sublattice $L'$ of $\ZZ^n$,
\item $L(D)$ a sublattice of $\ZZ^n$ and
\item $L'\oplus_{\ZZ} L(D) = \ZZ^n$.
\end{itemize}

(ii) A {\em Delaunay simplex set} $S$ to be a Delaunay polyhedron with $\vert S\vert = n+1$.

(iii) A {\em repartitioning set} $R$ to be a Delaunay polyhedron with $\vert S\vert = n+2$.

\end{definition}
The {\em isotropy lattice} $L(D)$ is uniquely determined by $D$ and its dimension is called the {\em degeneracy rank} of $D$ denoted $\degrk(D)$.
Note that $D$ is the vertex set of a convex body only when $L(D)=0$.
Also Delaunay polyhedra is full-dimensional, i.e. the smallest affine saturated lattice containing $D$ is $\ZZ^n$ itself.

The set $P_{L'}(D)$ is included in $L'$ and depends on $L'$.
For any two choices $L'_1$ and $L'_2$ there exist a bijective affine map
$\phi:L'_1\rightarrow L'_2$ with $\phi(P_{L'_1}(D)) = P_{L'_2}(D)$.
When we consider properties that do not depend on the integral representation,
we drop the lattice and write $P(D)$.

For $f\in Erdahl(n)$ we write
\begin{equation*}
Z(f)=\left\{x \in \ZZ^n\mbox{~such~that~} f(x)=0\right\}.
\end{equation*}

In classical geometry of numbers the essential tool is the quadratic form $Q$ instead of the quadratic function.
The following establish a direct link between both:
\begin{definition}
For a Delaunay simplex set $S \subset \ZZ^n$ and $Q\in S^n$ there exists a unique function $f\in E_2(n)$ such that:
\begin{itemize}
\item $f(x) = 0$ for $x\in S$
\item and $Q = \Quad(f)$.
\end{itemize}
This function is denoted $f_{S,Q}$ and depends linearly on $Q$.
\end{definition}

The key reason for using Delaunay polyhedra is the following theorem:
\begin{theorem}

(i) If $D$ is a Delaunay polyhedron then there exist a function $f\in Erdahl(n)$ such that $D=Z(f)$.

(ii) If $f\in Erdahl(n)$ then either $Z(f)$ is empty or there exist a $k$-dimensional saturated affine lattice $L\subset \ZZ^n$ such that $Z(f)$ is a Delaunay polyhedron of $L$.
\end{theorem}
\proof (i) Let us take a Delaunay polyhedron $D = P_{L'}(D) + L(D)$ with $P_{L'}(D)$ being a Delaunay polyhedron of
a lattice $L'$ with $L' \oplus_{\ZZ} L(D) = \ZZ^n$.
Let us denote by $S(c,r)$ the sphere around $P_{L'}(D)$ and write $L' = \ZZ v_1+\dots+\ZZ v_k$. The function
\begin{equation*}
\begin{array}{rcl}
f':\ZZ^n          & \rightarrow & \RR\\
x=(x_1,\dots,x_n) & \mapsto     & \left\Vert \sum_{i=1}^n x_i v_i -c \right\Vert^2-r^2
\end{array}
\end{equation*}
belongs to $Erdahl(k)$.
More precisely, $f'(x)=0$ if and only if $\sum_{i=1}^{k} x_i v_i \in P_{L'}(D)$.
For $x\in \ZZ^n$ we write $x= x_1 + z$ with $x_1\in L'$ and $z\in L(D)$ and write $f(x)=f'(x_1)$.
It is easy to prove that $f\in Erdahl(n)$ and $Z(f) = D$.

(ii) This is \cite[Corollary 2.5]{Erdahl1992}. \qed

We define the {\em rational closure} $S^{n}_{rat, \geq 0}$ to be the set of positive semidefinite
forms whose kernel is defined by rational equalities.

\begin{corollary}
If $f\in Erdahl(n)$ is such that $Z(f)$ is a Delaunay polyhedron then $\Quad(f)\in S^{n}_{rat, \geq 0}$.
\end{corollary}
\proof Let us write $D=Z(f)$ and take a lattice $L' \subset \ZZ^n$ with $L'\oplus_{\ZZ} L(D) = \ZZ^n$.
We write any $x\in \ZZ^n$ as $x = x_1 + z$ with $x_1\in L'$ and $z\in L(D)$.
There is a quadratic function $f_1$ on $L'$ such that $f(x) = f_1(x_1)$ and $Z(f_1) = P_{L'}(D)$.
Thus $\Quad(f_1)$ is positive definite and since $L(D)$ is an integral lattice the matrix $\Quad(f)$ belongs
to $S^{n}_{rat, \geq 0}$. \qed

Given a set $V\subset \ZZ^n$ we will need to be able to test if it is a Delaunay polyhedron or not.
Algorithm \ref{TestRealizability} does this iteratively for a finite point set by solving larger and
larger linear programs until conclusion is reached. 
The algorithm can be easily adapted to the case of a point set of the form $R + L$ with $L$ a lattice
and $R$ finite.
The corresponding algorithm for perfect form is given in \cite[Algorithm 1]{MinkowskianSublattices}.

\begin{algorithm}\label{TestRealizability}
\KwData{A finite set $V\subset \ZZ^n$ of full affine rank}
\KwResult{A quadratic function $f\in Erdahl_{>0}(n)$ such that $Z(f)=V$ if it exists and {\bf false} otherwise}
$S_{vert} \leftarrow V$.\\
$S_{vect} \leftarrow \emptyset$.\\
\Repeat{The set $\{ev_v\mbox{~for~}\in S_{vert}\}$ has rank $(n+1)(n+2)/2$}{
  $v \leftarrow \mbox{random~element~of~}\ZZ^n$\\
  $S_{vert} \leftarrow S_{vert}\cup \{v\}$
}
\While{no solution has been reached}{
  Form the linear program
  \begin{equation*}
  \begin{array}{rl}
  \mbox{minimize}_{f\in E_2(n)}   & Tr(\Quad(f))\\
  \mbox{subject~to} & f(v)=0     \mbox{~for~}v\in V\\
                    & f(v)\geq 1 \mbox{~for~}v\in S_{vert} - V\\
                    & \Quad(f)[v] \geq 1 \mbox{~for~}v\in S_{vect}\\
  \end{array}
  \end{equation*}
  \If{The linear program is infeasible}{
    return {\bf false}
  }
  $f_0 \leftarrow $ a rational optimal solution.\\
  \If{$f_0\in Erdahl(n)$ and $Z(f_0) = V$}{
    return $f_0$
  }
  \If{$\Quad(f_0)\notin S^n_{>0}$}{
    Find a vector $v\in \ZZ^n$ with $Q(f_0)[v]\leq 0$\\
    $S_{vect} \leftarrow S_{vect}\cup \{v\}$
  }
  \If{$\Quad(f_0)\in S^n_{>0}$ and $f_0\notin Erdahl(n)$}{
    Find a vector $v\in \ZZ^n$ with $f_0(v) < 0$\\
    $S_{vert} \leftarrow S_{vert}\cup \{v\}$
  }
  \If{$f_0\in Erdahl(n)$ and $Z(f_0)\not= V$}{
    Find a vector $v\in Z(f_0) - V$.\\
    $S_{vert} \leftarrow S_{vert}\cup \{v\}$
  }
}
\caption{Testing Delaunay realizability of a finite set of points}
\end{algorithm}

If $D$ is a $n$-dimensional Delaunay polyhedron, then we define
\begin{equation*}
\Aut(D)=\{\phi\in \AGL_n(\ZZ)\mbox{~:~}\phi(D)=D\}.
\end{equation*}
When using Algorithm \ref{TestRealizability} it is best to impose that the sought function $f$ is invariant under $\Aut(D)$ since it simplifies the search and a Delaunay polyhedron admits an invariant function (see Corollary \ref{InvFct} below).

Before stating our result on the description of $\Aut(D)$, we remind on the notion of semidirect product.
Given a group $G$, we call $G$ a semidirect product and write $G=N\times H$ if $N$ is a normal subgroup
of $G$, $H$ a subgroup $G=NH$ and $N\cap H=\{e\}$.

\begin{theorem}\label{AutomAut_D}
If $D$ is a $n$-dimensional Delaunay polyhedron of degeneracy degree $d$ then
we have the isomorphism
\begin{equation*}
\Aut(D)=\left\langle (\ZZ^d)^{1+n-d} \rtimes \GL_d(\ZZ)\right\rangle \rtimes \Aut(P(D)).
\end{equation*}

\end{theorem}
\proof Let us take a basis $v=(v_i)_{1\leq i\leq d}$ of $L(D)$.
Any automorphism of $D$ will send $v$ to another basis of $L(D)$ and this determines a component $\GL_d(\ZZ)$ of the automorphism group.
Let us write $L'$ for an affine sublattice of $\ZZ^n$ such that $L'\oplus_{\ZZ} L(D) = \ZZ^n$. Such affine sublattice are determined by $1+n-d$ vectors in $L(D)$ and this determines the component $(\ZZ^d)^{1+n-d}$ of the automorphism group.
The last component comes from the fact that the automorphisms have to preserve the polytope $\conv(P(D))$. \qed

We denote by $\Aff(D)$ the normal subgroup $(\ZZ^d)^{1+n-d} \rtimes \GL_d(\ZZ)$ given in the above theorem.

\begin{corollary}\label{InvFct}
A Delaunay polyhedron $D$ admits a function $f\in Erdahl(n)$ with $Z(f)=D$ that is invariant under $\Aut(D)$.
\end{corollary}
\proof Let us write $L' \oplus_{\ZZ} L(D) = \ZZ^n$. A vector $x\in \ZZ^n$ is decomposed as $x=x_1 + z$
with $x_1\in L'$ and $z\in L(D)$.
Any function $f$ with $Z(f) = D$ must have $\Quad(f)[z] = 0$ and $\Lin(f) \cdot z=0$.
Therefore, there exist a function $f_1$ on $L'$ such that $f(x) = f_1(x_1)$. So, $f$ is invariant under $\Aff(D)$.
The group acting on $f_1$ is $\Aut(P_{L'}(D))$, which is finite.
The function
\begin{equation*}
f'_1(x_1) = \sum_{u\in \Aut(P_{L'}(D))} f(u(x_1))
\end{equation*}
is $\Aut(P_{L'}(D))$ invariant.
Thus we get a function $f'(x) = f'_1(x_1)$ invariant under $\Aut(D)$. \qed

If $D$ and $D'$ are two Delaunay polyhedra such that $D\subset D'$ then we define the stabilizer group
\begin{equation*}
\begin{array}{rcl}
  \Stab(D, D') &=& \left\{\phi \in \AGL_n(\ZZ)\mbox{~:~}\phi(D)=D\mbox{~and~}\phi(D')=D'\right\}\\
               &=& \Aut(D) \cap \Aut(D').
\end{array}
\end{equation*}

We have the following results
\begin{theorem}\label{InclusionAff}
Let $D$ and $D'$ be two Delaunay polyhedra satisfying $D\subset D'$.
Then:
\begin{enumerate}
\item[(i)] we have $L(D)\subset L(D')$ and $\Aff(D) \subset \Aff(D')$.
\item[(ii)] There exist a finite group $G_1\subset \Aut(P(D))$ such that
\begin{equation*}
\Stab(D, D') = \Aff(D) \rtimes G_1
\end{equation*}
in particular $\Aff(D)$ is a finite index subgroup of $\Stab(D, D')$
\end{enumerate}

\end{theorem}
\proof We have $L(D)\subset L(D')$.
Let us take a $\ZZ$-basis $\{e_1, \dots, e_d\}$ of $L(D)$ and complement it by adding $\{e_{d+1}, \dots, e_{d'}\}$ to a basis of $L(D')$.
We can find $f_1, \dots, f_{n-d'}$ so that the $e_i$ and $f_j$ form a basis of $\ZZ^n$.
The group $\Aff(D')$ is generated by translations along $f_1, \dots, f_{n-d'}$ and $\GL_{d'}(\ZZ)$.
The group generated by translations along $\{e_{d+1}, \dots, e_{d'}\}$ and $\GL_d(\ZZ)$ directly embeds into $\GL_{d'}(\ZZ)$ and this determine the group inclusion. So (i) holds.

By (i) we have the inclusion $\Aff(D) \subset \Stab(D,D')$.
Let us choose a lattice $L'$ with $L'\oplus_{\ZZ} L(D) = \ZZ^n$. 
Then if $u\in \Stab(D,D')\subset \Aut(D)$ we can find an unique element
$n\in \Aff(D)$ such that $u n^{-1}$ stabilizes $L'$.
Thus $u n^{-1}$ belongs to $\Aut(P_{L'}(D))$. The image defines the group $G_1$.
the finite index property follows. \qed

\section{Facial structure of the Erdahl cone}\label{FacialErdahlCone}
The standard scalar product on $S^n$ is $\langle A, B\rangle = Tr(AB)$.
We equip $Erdahl(n)$ with the inner product
\begin{equation*}
(f, g) =  \Cst(f) \Cst(g) + 2 \Lin(f) \cdot \Lin(g)  + \langle \Quad(f), \Quad(g) \rangle
\end{equation*}
and for each $x\in \ZZ^n$ we define the evaluation function $ev_x$ by
$ev_x(y) = (1 + x\cdot y)^2$ such that
\begin{equation*}
(f, ev_x) =  f(x).
\end{equation*}

A \textit{convex cone} $\MC$ is defined as a set invariant
under addition and multiplication by positive scalars.
$\MC$ is called \textit{full-dimensional} if the only
vector space containing it is $\RR^m$.
$\MC$ is called \textit{pointed} if no linear subspace of
positive dimension is contained in it.
Let $\MC$ be a full-dimensional pointed convex polyhedral cone in $\RR^m$.
Given $f \in (\RR^{m})^{*}$, the inequality $f(x)\geq 0$
is said to be \textit{valid} for $\MC$ if it holds for all $x\in \MC$.
A \textit{face} of $\MC$ is a pointed polyhedral cone
$\{ x \in \MC \mbox{~}:\mbox{~} f(x) = 0 \}$,
where $f(x) \geq 0$ is a valid inequality.

A face of dimension $1$ is called an \textit{extreme ray} of $\MC$;
a face of dimension $m-1$ is called a \textit{facet} of $\MC$.
The set of faces of $\MC$ forms a partially ordered set under inclusion.
We write $F\lhd G$ if $F\subset G$ and $\dim F=\dim G-1$.
Two extreme rays of $\MC$ are said to be \textit{adjacent}
if they generate a two-dimensional face of $\MC$.
Two facets of ${\mathcal C}$ are said to be \textit{adjacent}
if their intersection has dimension $m - 2$.
Any $(m-2)$-dimensional face of $\MC$ is called a \textit{ridge}
and it is the intersection of exactly two facets of $\MC$.

By the Farkas-Minkowski-Weyl Theorem (see e.g.\ \cite[Corollary~7.1a]{schrijver}) a convex cone $\MC$ is \textit{polyhedral} if and only it is defined either 
by a finite set of \textit{generators} $\{v_1,\ldots, v_N\} \subseteq \RR^m$
or by a finite set of linear functionals $\{f_1, \ldots, f_M\}\subseteq (\RR^m)^*$:
\[
\MC =
\bigg\{\sum_{i=1}^N \lambda_i v_i \mbox{~}:\mbox{~} \lambda_i\geq 0\bigg\} =
\bigg\{x\in \RR^m \mbox{~}:\mbox{~} f_i(x)\geq 0\bigg\}.
\]

Every minimal set of generators $\{v_1, \ldots, v_{N'}\}$ defining
a polyhedral cone $\MC$ has the property 
\[\{\RR_+ v_1,\ldots, \RR_+ v_{N'}\}=\{e \mbox{~}:\mbox{~} \mbox{$e$ extreme ray of $\MC$}\}.\]
Every minimal set of linear functionals $\{f_1, \ldots, f_{M'}\}$
defining $\MC$ has the property that $\{F_1, \ldots, F_{M'}\}$
with $F_i=\{x\in \MC \mbox{~}:\mbox{~} f_i(x)=0\}$ is the set of
facets of $\MC$.
The problem of transforming a minimal set of generators into
a minimal set of linear functionals (or vice versa) is called
the \textit{dual description problem}.

In our work, we have to deal with Delaunay polyhedra with an infinite number
of vertices and we cannot apply the Farkas-Minkowski-Weyl theorem to them nor
of course existing dual-description software \cite{lrs,cdd}.

\begin{definition}
Let $D\subset \ZZ^n$ be a Delaunay polyhedron

(i) We define the vector space
\begin{equation*}
Space(D)=\{f\in E_2(n)\mbox{~such~that~} f(x)=0\mbox{~for~}x\in D\}.
\end{equation*}

(ii) The dimension of $Space(D)$ is called the {\em perfection rank} $\rankperf(D)$
and $D$ is {\em perfect} if $\rankperf(D)=1$.
\end{definition}

\begin{proposition}
Let $D\subset \ZZ^n$ be a Delaunay polyhedron 

(i) $Space(D) \cap Erdahl(n)$ is a face of $Erdahl(n)$ of dimension $\dim Space(D)$.

(ii) If $D$ is perfect then $Space(D) \cap Erdahl(n)$ is an extreme ray of $Erdahl(n)$.
\end{proposition}
\proof Let $p=\dim Space(D)$ and $g_1$, \dots, $g_p$ be a basis of $Space(D)$.
Since $D$ is a Delaunay polyhedron there exist a function $f\in Erdahl(n)$ such that $D=Z(f)$. 
For each $1\leq i\leq p$ there exist $\lambda_i > 0$ such that $\lambda_i f + g_i\in Erdahl(n)$
and so (i) follows.
(ii) follows directly from (i). \qed

For a perfect Delaunay polyhedron $D$ we denote $f_D$ a generator of the extreme ray $Space(D)\cap Erdahl(n)$.
\begin{theorem}\label{ClassifKindExtRay}(\cite[Theorem 2.1]{Erdahl1992})
The generators of extreme rays of $Erdahl(n)$ are:
\begin{enumerate}
\item The constant function $f=1$
\item The functions of the form $f_{a,\beta}(x)=(a_1 x_1 + \dots a_n x_n + \beta)^2$ with $(a_1, \dots, a_n)$ not colinear to an integral vector.
\item The functions of the form $f_D$ with $D$ a perfect Delaunay polyhedron.
\end{enumerate}
\end{theorem}

This theorem indicates that the structure of the extreme rays of $Erdahl(n)$ is more complicated
than for a polytope.
Since we are interested only in the third class of extreme rays, some reduction will be necessary
and it turns out that we can work out everything with Delaunay polyhedra.

In this paper we will work with both spaces of functions in $E_2(n)$ and with point sets of
Delaunay polyhedra:
\begin{definition}
Given a Delaunay polyhedron $D$ we define

(i) The cone of admissible functions is defined as
\begin{equation*}
Erdahl_{supp}(D)=\left\{g\in E_2(n)\;:\; g(x)\geq 0 \mbox{~for~all~} x\in D\right\}.
\end{equation*}

(ii) The cone of evaluation functions is defined as
\begin{equation*}
Erdahl_{supp}^*(D)= \left\{ \sum_{x\in D} \lambda_x ev_x \mbox{~with~} \lambda_x\geq 0\right\}.
\end{equation*}
\end{definition}

\begin{theorem}\label{BasicResultErdahlCone}
Let $D$ be a Delaunay polyhedron of perfection rank $r$.

(i) $Erdahl_{supp}(D)$ is the product of a pointed convex cone $C_D$ with $Space(D)$.\\
The dual of $C_D$ is $Erdahl_{supp}^*(D)$.

(ii) Any Delaunay polyhedron $D'\subset D$ of perfection rank $r+1$ gives a facet of $Erdahl_{supp}^*(D)$.

(iii) If $L(D) = 0$ then facets of $Erdahl_{supp}^*(D)$ correspond to Delaunay polyhedra $D'\subset D$ of perfection rank $r+1$.
\end{theorem}
\proof (i) If $f$ and $-f$ both belong to $Erdahl_{supp}(D)$ then $f(x)=0$ for $x\in D$ and so $f\in Space(D)$.
So, $Erdahl_{supp}(D)$ is the sum of $Space(D)$ and a closed convex cone.
The duality result follows from \cite[Part IV.5]{barvinok}) for closed full-dimensional convex cones.

(ii) Let $f\in Erdahl(n)$ such that $Z(f) = D'$.
Thus we have $f(x)=0$ on $D'$ and $f(x) > 0$ on $D - D'$ and so $f$ defines a facet of $Erdahl_{supp}^*(D)$.

(iii) If $L(D) = 0$ then for any facet of $Erdahl_{supp}^*(D)$ we can find a set $D'\subset D$ and a function $f\in E_2(n)$ with $f(x)=0$ for $x\in D'$ and $f(x)>0$ for $x\in D - D'$. There exist a function $g\in Erdahl(n)$ such that $D=Z(g)$. Then we can find $\lambda > 0$ such that $f + \lambda g\in Erdahl(n)$. Then we have $D'=Z(f + (\lambda+1) g)$ and so $D'$ is a Delaunay polyhedron. \qed

Theorem \cite[Theorem 2.1]{Erdahl1992} shows that if $L(D)\not= 0$ there are other facets of $Erdahl_{supp}^*(D)$
than the ones from Delaunay polyhedra.

\begin{proposition}\label{LowerBoundNrVertices}
If $D$ is a $n$-dimensional perfect Delaunay polyhedron with degeneracy degree $d$ then $P(D)$ has at least
\begin{equation*}
{n-d + 2 \choose 2} -1 
\end{equation*}
points.
\end{proposition}
\proof The isotropy lattice $L(D)$ has dimension $d$ and we can choose a complement lattice $L'$
such that $L' \oplus_{\ZZ} L(D) = \ZZ^n$.
If $f$ is a quadratic function with $Z(f) = D$ then $f$ is determined by its restriction to $L'$.
Hence $f$ belongs to a vector space of dimension ${n-d + 2 \choose 2}$ and this gives the minimal
number of determining inequalities. \qed

A surprising relation has been found between the Erdahl cone of the hypercube
$\{0,1\}^n$ and the cut polytope, which is classic polytope of Combinatorial Optimization \cite{DL}:
Write $N=\{1,\dots, n\}$; if $S\subset N$,
then the cut metric $\delta_S$ on $N$ is defined as follows:
\begin{equation*}
\delta_S(i,j)=\left\{\begin{array}{rl}
1 & \mbox{ if } |\{ij\}\cap S|=1,\\
0 & \mbox{ otherwise}.
\end{array}\right.
\end{equation*}
We have $\delta_S=\delta_{N-S}$ and the {\em cut polytope} $\CUTP_{n}$
is defined as the convex hull of the cut metrics $\delta_S$.
The cone defined by the cut polytope is defined as
\begin{equation*}
\CCUTP_n = \left\{ \sum_{S\subset \{1,\dots,n\}} \lambda_S (1, \delta_S) \mbox{~with~} \lambda_S \geq 0\right\}.
\end{equation*}
The facets of the cone $\CCUTP_n$ are in one-to-one correspondence with the facets of the polytope $\CUTP_n$.

\begin{theorem}\label{CUTn_and_ErdahlHn}
The polyhedral cone $Erdahl_{supp}^*(\{0,1\}^n)$ is linearly equivalent to $\CCUTP_{n+1}$.
\end{theorem}
\proof The hypercube $[0,1]^n$ is defined as the convex hull of $2^n$ vectors
$v=(v_1, \dots, v_n)$ with $v_i\in \{0,1\}$.
For every such vector the evaluation function is
\begin{equation*}
ev_v(x)=1 + 2(v\cdot x) + (v\cdot x)^2.
\end{equation*}
Thus we can associate to $ev_v$ the vector
\begin{equation*}
(1, (v_i)_{1\leq i\leq n}, (v_iv_j)_{1\leq i\leq j\leq n}).
\end{equation*}
Since $v_i^2=v_i$ for $v_i\in \{0,1\}$ this vector family is linearly equivalent to
\begin{equation*}
(1, (v_i)_{1\leq i\leq n}, (v_i v_j)_{1\leq i < j\leq n}) \mbox{~for~}v\in \{0,1\}^n .
\end{equation*}
To the same $v$ we can associate the vector $\overline{v}=(0,v_1,\dots,v_n)$
and the set $S=\{i\in \{0,\dots, n\}, \mbox{~s.t.~} \overline{v}_i=1\} \subset \{0,\dots,n\}$.
The cut metric $\delta_S$ on $\{0,\dots, n\}$ is characterized by
\begin{equation*}
(\delta_S(i,j) )_{0\leq i<j\leq n}
\end{equation*}
and since $\delta_S(i,j)=(\overline{v}_i - \overline{v}_j)^2$ the family of pairs $(1, \delta_S)$ is linearly equivalent to the family of evaluation map $ev_v$. \qed

It is interesting to note that the symmetry group of the hypercube $[0,1]^n$ is
of size $2^n n!$ but that the symmetry group of the cut polytope $\CUTP_{n+1}$
is of size $2^n (n+1)!$ for $n\not=3$ \cite{Relatives,DL}.

\section{The Delaunay polyhedra retract}\label{GeometryErdahlCone}

For $f\in Erdahl(n)$, we define $Vect\,Z(f)$ to be the vector space
spanned by difference of elements of $Z(f)$.
We define $V_f=Vect\, Z(f) + \Ker\, \Quad(f)$.
For a given $f\in Erdahl(n)$ we call {\em proper pair} a pair $(g,h)\in E_2(n)^2$ such that $g\in Erdahl(n)$, $h(x)\geq 0$ for $x\in \RR^n$ and $f=g + h$.

\begin{lemma}\label{KernelsOfH}
Let $f\in Erdahl(n)$.

(i) For a proper pair $(g,h)$, $Z(f)\subset Z(g)$ and $\Ker\,\Quad(h)\subset V_f$.

(ii) There exist a proper pair $(g,h)$ with $\Ker\,\Quad(h) = V_f$

\end{lemma}
\proof If $x\in Z(f)$ then one sees that necessarily $h(x)=g(x)=0$.
So, $h(x)=0$ for $x\in \conv(Z(f))$ which implies $Vect\, Z(f)\subset Ker\, \Quad(h)$.
Also, it is clear that for any vector $v\in \Ker\,\Quad(f)$ we have $\Quad(g)[v]=\Quad(h)[v]=0$. Hence, (i) holds.

Let us denote by $\ZZ\Ker\,\Quad(f)$ the smallest subspace of $\RR^n$ having an integral basis containing $\Ker\,\Quad(f)$. 
By the Decomposition Lemma 3.1 in \cite{Erdahl1992} there exist a $g\in Erdahl(n)$ with $\Ker\,\Quad(g) = \ZZ\Ker\,\Quad(f)$ and a positive semidefinite form $Q_1$ such that
\begin{equation*}
f(x) = g(x) + Q_1[x].
\end{equation*}
Let us denote by $V$, respectively $W$, an integral supplement of $\Ker\, \Quad(f)$ in $\ZZ\Ker\,\Quad(f)$, respectively $\ZZ\Ker\,\Quad(f)$ in $\ZZ^n$.
Denote by $\phi_1$, \dots, $\phi_m$ some affine functions on $\ZZ^n$ such that $\phi_i(\ZZ\Ker\,\Quad(f))=0$ and 
\begin{equation*}
\{x\in W\mbox{~s.t.~} \phi_1(x)=\dots=\phi_m(x)=0\} = Vect\,Z(f) \cap W.
\end{equation*}
Then for $\epsilon>0$ small enough the function $g_1 - \epsilon\sum_{i=1}^m \phi_i(x)^2$ is still in $Erdahl(n)$.
So, one gets that the pair $(f-h,h)$ with $h(x) = Q_1[x] + \epsilon\sum_{i=1}^m \phi_i(x)^2$ is proper and (ii) is true. \qed

Let us call $W$ an integral supplement of $V_f$. Denote by $\Quad(f)\vert_{W}$ the quadratic form $\Quad(f)$ restricted to $W$.
A proper pair $(g,h)$ is called {\em extremal} if $\det\,\Quad(h)|_{W}$ is maximal among all proper pairs.
Lemma \ref{KernelsOfH}.(i) implies that the notion of being extremal is
independent of the chosen subspace $W$, while Lemma \ref{KernelsOfH}.(ii)
implies that there is at least one form of non-zero determinant.

\begin{theorem}\label{UniqueDecomposition}
Let $f\in Erdahl(n)$.

(i) An extremal proper pair has $Z(g)$ being a Delaunay polyhedron.

(ii) There exist a unique extremal proper pair $(g,h)$.

\end{theorem}
\proof Let us take an integral supplement $W$ as above and suppose to avoid trivialities that $W\not=\emptyset$.
So, by restricting $f$ to $W$, we can assume that $\Quad(h)$ is positive definite.
\begin{equation*}
f(x) = g(x) + h(x)
\end{equation*}
with $h(x)\geq 0$ for $x\in \RR^n$. This condition on $h$ is equivalent to 
\begin{equation*}
A_h=\left(\begin{array}{cc}
\Cst(h)  & \Lin(h)^T\\
\Lin(h)  & \Quad(h)
\end{array}\right)
\end{equation*}
being positive definite.
Hence we consider the following semidefinite programming problem: find the $A_h\in S^{n+1}_{\geq 0}$ maximizing $\det\,\Quad(h)$ and satisfying for all $x\in \ZZ^n$
\begin{equation*}
f(x) \geq A_h[(1,x)] = h(x).
\end{equation*}
We also write $g=f-h$.

Suppose that $Z(g)$ is not a Delaunay polyhedron and that $\Quad(h)$ is positive definite.
Then $Z(g)$ does not generates $\RR^n$ as an affine space and so there
exists an affine function $\phi$ such that $\phi( Z(g) )=0$.
Then there exist $\alpha > 0$ such that the pair $(f-h', h')$ with $h'= h + \alpha \phi^2$ is still proper.
Since $\det\,\Quad(h') > \det\,\Quad(h)$ the pair is not extremal and (i) holds.

Let us take $N=n(n+1)/2$ points  $v_i\in \ZZ^n$ such that the family $\left\{(1,v_i)(1,v_i)^T\right\}_{1\leq i\leq N}$ is of full rank. The inequalities $g(v_i) \geq A[(1,v_i)] \geq 0$ implies that all coefficients of $A$ are bounded. 
Thus the problem is actually to minimize the convex function $h\mapsto -\log\det\,\Quad(h)$ over a compact convex set, hence existence follows.

Since $-\log\, \det$ is a strictly convex function, we know that if we have two optimal solutions $h_1$ and $h_2$ then $\Quad(h_1) = \Quad(h_2)$.
Let us denote by $D_1=Z(g_1)$ and $D_2=Z(g_2)$ the corresponding Delaunay polyhedra.
The function $h_{mid} =(h_1 + h_2)/2$ is also an optimal solution of the problem.
We have $Z(g_{mid}) = D_1\cap D_2$.
The set $D_1\cap D_2$ is necessarily a Delaunay polyhedron since otherwise we could still increase the determinant by the above construction and this would contradict the optimality.
But if $Z(f)$ is a Delaunay polyhedron then the terms $\Cst(f)$ and $\Lin(f)$ are determined by $\Quad(f)$.
So, one gets $h_1=h_2$ and the uniqueness is proved on the restriction to $W$.
But Lemma \ref{KernelsOfH}(i) implies that once $\Quad(h)$ is known on $W$ then it is known on $\ZZ^n$.
By the condition $h\geq 0$ the linear part is known as well. \qed

Note that in the above determinant maximization problem a finite set of
inequalities suffices to determine the optimal solution.
This follows from the fact that since we are maximizing the determinant
we can assume that the lowest eigenvalue of $\Quad(h)$ is bounded away
from $0$, i.e. that there exist $c>0$ such that $\Quad(h) \geq c I_n$.

For $f\in Erdahl(n)$, we write $proj(f)=g$ and $proj'(f)=h$ with $(g,h)$ the unique extremal pair associated to $f$.
From the unicity of extremal pairs we also get that $proj$ and $proj'$ commute with the action of $\AGL_n(\ZZ)$.

\begin{conjecture}
The function $proj$ is continuous.
\end{conjecture}
Let us define $Erdahl_{dp}(n)$ to be the set of $f\in Erdahl(n)$ such that $Z(f)$ is a Delaunay polyhedron.
The above conjecture if true implies that the set $Erdahl_{dp}(n)$ is simply connected and this could be of interest for topological applications.
However, we were not able to prove the conjecture and instead we prove in later sections connectedness results which are sufficient for our purposes.

\section{Relation with $L$-types theory}\label{SectionLtype}

In this section we reframe classical $L$-type theory from \cite{VoronoiII}
(see also \cite{bookschurmann} for a modern account) in term of Erdahl cone
and state several key lemmas.

\begin{definition}
  Let $Q \in S^n_{rat,\geq 0}$. The Delaunay polyhedra tessellation ${\mathcal DPT}(Q)$ defined by
  $Q$ is the set of Delaunay polyhedra $D$ such that there exist a $f\in Erdahl(n)$ with
  \begin{itemize}
  \item $Z(f) = D$
  \item and $\Quad(f) = Q$.
  \end{itemize}
\end{definition}
If $Q$ is positive definite then the Delaunay polyhedra tessellation is the classical
Delaunay polytope tessellation, i.e. all Delaunay polyhedra occurring are actually vertex
sets of Delaunay polytopes.
The number of translation classes of Delaunay polyhedra is always finite.
These Delaunay polyhedra tessellations were considered in \cite[Section 2.2]{EquivariantLtypeDSV}.
Efficient algorithm for the enumeration of Delaunay polytope tessellations are given
in \cite{ComplexityVoronoiDSV}.

From this one can define the $L$-type which are parameter spaces of Delaunay polytope tessellations:

\begin{definition}
Let us take a Delaunay polyhedra tessellation ${\mathcal T}$. Then the $L$-type $LT({\mathcal T})$
is defined as the closure of the set of quadratic forms $Q$ such that ${\mathcal DPT}(Q) = {\mathcal T}$.
It is well known (see \cite{bookschurmann} and \cite{VoronoiII} for proofs) that $L$-types are polyhedral
cones.
\end{definition}
A $L$-type is called {\em primitive} if it is of maximal dimension, 
this is equivalent to say that all its Delaunay polyhedra are Delaunay simplex sets.

The set of all $L$-types for all possible Delaunay tessellations defines
a tessellation of the cone $S^{n}_{rat, \geq 0}$.

Given two Delaunay polyhedra tessellation ${\mathcal T}$ and ${\mathcal T}'$ we say that
{\em ${\mathcal T}'$ is a refinement of ${\mathcal T}$} if every Delaunay polyhedron of ${\mathcal T}'$
is included in a single Delaunay polyhedron of ${\mathcal T}$.
${\mathcal T}'$ is a {\em simplicial refinement} if all its Delaunay polyhedra are Delaunay simplex sets.

\begin{proposition}\label{ExistenceSimplicialRefinement}
Any Delaunay polyhedra tessellation ${\mathcal T}$ admits at least one simplicial refinement.
\end{proposition}
\proof Let us denote by $L({\mathcal T})$ the space $L(D)$ of the Delaunay polyhedra $D$ occurring
in the tessellation and by $Q \in S^n_{rat, \geq 0}$ the form realizing it.
Let us take a lattice $L'$ such that $L' \oplus_{\ZZ} L({\mathcal T}) = \ZZ^n$.
We write $x\in \ZZ^n$ as $x=x_1 + z$ with $x_1\in L'$ and $z\in L({\mathcal T})$.
$D$ is a Delaunay polyhedron for the quadratic function $f\in Erdahl(n)$.
Necessarily $f$ is of the form $f(x) = f_1(x_1)$ with $f_1$ a quadratic function
on $L'$.
Let us denote $(D_i)_{i\in I}$ the Delaunay polyhedra occurring in the tessellation.
Let us take a basis $w_1, \dots, w_m$ of $L({\mathcal T})$ and define linear forms $\phi_i$ on $\ZZ^n$
such that $\phi_i(w_j) = \delta_{ij}$ and $\phi_i(L') = 0$. The quadratic form
\begin{equation*}
Q'[x] = Q[x] + \sum_{i=1}^m  \left(\phi_i(x)\right)^2
\end{equation*}
is positive definite and the Delaunay polyhedra tessellation ${\mathcal DT}$
corresponding to $Q'$ is formed by the Delaunay polyhedra
\begin{equation*}
D_i + \sum_{k=1}^m \{a_k, a_k+1\} w_k \mbox{~with~}1\leq i\leq r \mbox{~and~} a_k\in \ZZ.
\end{equation*}
In particular ${\mathcal T}_2$ is a refinement of ${\mathcal T}$.

Since the $L$-type domain form a tiling of $S^n_{rat, \geq 0}$ the form $Q'$ belongs
to at least one primitive $L$-type $LT({\mathcal T}_2)$.
This $L$-type defines a Delaunay polyhedra tessellation by simplices
which is a refinement of ${\mathcal T}_2$ and so of ${\mathcal T}$. \qed

Let us take a primitive $L$-type ${\mathcal T}$.
Any facet $F$ of ${\mathcal T}$ is determined by a pair of Delaunay simplex sets $S_1$
and $S_2$ in the Delaunay tessellation that determine a repartitioning set.
We say that two facet defining repartitioning sets are in the same class if they
define the same facet $F$ of ${\mathcal T}$.
If $R$ is a repartitioning set then $\conv(R)$ admits exactly two triangulations
(see \cite[Section 4.3.2]{bookschurmann}).
One says that two primitive $L$-types are adjacent if their intersection
is a codimension $1$ face in the cone $S^n_{>0}$.
When we move from one $L$-type to another $L$-type, the Delaunay tessellation
is changed and this is done combinatorially by the repartitioning sets.
That is some Delaunay simplex sets are merged into repartitioning sets
and the triangulation is changed to the other triangulation, thus yielding
another $L$-type.

Given a Delaunay polyhedron $D$ a Delaunay polyhedra tessellation ${\mathcal T}$
is called {\em $D$-proper} if $D$ is the union of the Delaunay polyhedra $D'$ contained in $D$.
We have following lemma:

\begin{lemma}\label{ConvexityDproper}
  Let $D$ be a Delaunay polyhedron.
  The graph formed by the primitive $L$-types whose corresponding Delaunay polyhedra tessellations
  is primitive and $D$-proper is connected.
\end{lemma}
\proof Let us consider a function $f_D\in Erdahl(n)$ such that $Z(f_D) = D$.
We can consider the triangulations induced by positive definite quadratic forms on $D$ itself.
This set is connected by the theory of regular triangulations
(see \cite{TriangulationBook} for an account).

Any triangulation ${\mathcal T}_{part}$ on $D$ induced by a positive definite quadratic form $Q$
can be extended to a triangulation ${\mathcal T}$ of $\ZZ^n$:
It suffices to replace $Q$ by $Q + \lambda \Quad(f_D)$ for $\lambda$ sufficiently large.
The reason is that $\Quad(f_D)$ will not change the Delaunay triangulation for Delaunay simplex sets
contained in $D$.

Now given a primitive $L$-type $LT$ whose Delaunay polyhedra tessellation ${\mathcal T}$
is $D$-proper, we denote by ${\mathcal S}$ its set of Delaunay simplex sets included in $D$.
We consider the following cone ${\mathcal C}({\mathcal S})$:
\begin{equation*}
{\mathcal C}({\mathcal S}) = \left\{ Q\in S^n_{rat,\geq 0} \mbox{~s.t.~} f_{S, Q} (x) \geq 0 \mbox{~for~} S\in {\mathcal S}\mbox{~and~} x\in \ZZ^n - S\right\}.
\end{equation*}
This cone is convex and is an union of primitive $L$-types. Thus this set of $L$-types is connected.
The connectedness follows by combining above results. \qed

Delaunay polytopes restricted to the set $D$.
When doing so, by the connectivity

\begin{theorem}\label{ConnectivitySimplices}
Let us take $D$ a Delaunay polyhedron and two Delaunay simplex sets
$S$ and $S'$ in $D$. Then there exist a sequence of Delaunay simplex sets
$\{S=S_0,S_1, \dots, S_m=S'\}$ with $S_i\subset D$ for $0\leq i\leq m$
such that $S_i\cup S_{i+1}$ is a repartitioning set for $0\leq i\leq m-1$.
\end{theorem}
\proof Let us take $f_D\in Erdahl(n)$ a function such that $Z(f_D)=D$.
Take $f_S$, $f_{S'}$ the corresponding functions for $S$ and $S'$.
Denote by ${\mathcal T}$ the Delaunay polyhedra tessellation defined by 
$\Quad(f_D)$, which obviously has $D$ as one of its component.

When we perturb $\Quad(f_D)$, we are changing the Delaunay tessellation.
However, if we take $\epsilon > 0$ small enough, we can ensure that the
Delaunay polyhedra tessellations 
${\mathcal DPT}(\Quad(f_D + \epsilon f_S))$ and ${\mathcal DPT}(\Quad(f_D + \epsilon f_S'))$
are $D$-proper.
By applying Proposition \ref{ExistenceSimplicialRefinement}
we can find simplicial refinement of those two tessellations which we name ${\mathcal TR}$ and
${\mathcal TR}'$ and are both $D$-proper.
We call $LT$ and $LT'$ the corresponding primitive $L$-types.

By Lemma \ref{ConvexityDproper} there exist a path between $LT$ and $LT'$ that uses
only $D$-proper $L$-types.
By following this path, we can change $S$ into another
Delaunay simplex set $S_2$ in ${\mathcal TR}'$.

Denote by $f_1$, \dots, $f_r$ the facets of $LT'$.
Every such facet corresponds to a family of repartitioning sets.
We say that two Delaunay simplex sets included in $D$ are adjacent if their union
is a repartitioning set which gives a facet of $LT'$.
The Delaunay polyhedron $D$ is a coarsening obtained by merging all simplices,
so the above defined graph is connected.
This means that we can find a path from $S_2$ to $S'$. \qed

\section{Relation with hypermetric theory}\label{SectionHypermetric}

We define the {\em volume} $\Vol(S)$ of a Delaunay simplex set $S$ to be $n! Vol(\conv(S))$
with $Vol$ the Euclidean volume.
This rescaled volume is an integer and satisfies $\Vol(S)\leq n!$.
The possible rescaled volumes $PossVol(n)$ are given in
Table \ref{TableRescaledVolume} for $n\leq 7$ and a super-exponential lower
bound on $\max\, PossVol(n)$ is proven in \cite{SuperExpVolume}.
The best known upper bound \cite[Proposition 14.2.4]{DL} is
\begin{equation}\label{UpperBound_PossVol}
\max\,PossVol(n) \leq n! \frac{2^n}{{2n\choose n}}.
\end{equation}

\begin{table}
\label{TableRescaledVolume}
\caption{Possible volume of lattice Delaunay simplices. See \cite{Ba} for the proof for $n\leq 6$ and Section \ref{DimSeven_simplices} for the proof for $n=7$.}
\begin{tabular}{|c|c|c|c|c|c|c|c|}
\hline
$n$                & $1$     & $2$     & $3$     & $4$     & $5$       & $6$           & $7$\\
\hline
$PossVol(n)$ & $\{1\}$ & $\{1\}$ & $\{1\}$ & $\{1\}$ & $\{1,2\}$ & $\{1,2,3\}$   & $\{1,2,3,4,5\}$\\
\hline
\end{tabular}
\end{table}

\begin{definition}
Let us take two $n$-dimensional Delaunay polyhedron $D$, $D'$ with $D\subset D'$.
We can define the {\em generalized hypermetric cone}
\begin{equation*}
\Hyp(D,D')=\left\lbrace f\in E_2(n) \mbox{~s.t~} f(x)=0\mbox{~if~}x\in D
\mbox{~and~}f(x)\geq 0\mbox{~if~}x\in D'\right\rbrace .
\end{equation*}
\end{definition}
We have the inclusion $\Hyp(D,D')\subset Erdahl_{supp}(D')$ and $\Hyp(D,D')$ is a priori defined by
an infinity of inequalities.

As a direct application we can express the $L$-type domains as intersection of generalized
hypermetric cones:
\begin{proposition}
  Let ${\mathcal T}$ be a Delaunay polyhedra tessellation. Then we have
\begin{equation*}
LT({\mathcal T}) = \cap_{D\in {\mathcal T} } \Quad \Hyp(D, \ZZ^n).
\end{equation*}
  
\end{proposition}

\begin{proposition}
The polyhedral cone $\Hyp(D,D')$ is polyhedral.
\end{proposition}
\proof Let us take a Delaunay simplex set $S=\{v_0, \dots, v_n\} \subset D$ which exists
by Proposition \ref{ExistenceSimplicialRefinement}.
If we prove the polyhedrality of $\Hyp(S,D')$ then $\Hyp(D,D')$ is polyhedral as well
since it is obtained from $\Hyp(S, D')$ by adding equalities $f(x)=0$ for $x\in D - S$.

Suppose that a $v\in D'$ defines a relevant inequality. Then there exist a function $f$ such that
$f(x)=0$ for $x\in S\cup \{v\}$ and $f(x) >0$ for $x\in D' - S \cup \{v\}$.
Since $D'$ is a Delaunay polyhedron, there exist a function $g$ such that $g(x)=0$ for $x\in D'$
and $g(x) > 0$ for $x\in \ZZ^n - D'$. Then we can find $\lambda>0$ such that $f(x) + \lambda g(x) > 0$
for $x\in \ZZ^n - S \cup \{v\}$.
As a consequence the polytope $\conv(S \cup \{v\})$ is a Delaunay polytope.
This implies that for any $i\in \{0,\dots, n\}$ the Delaunay simplex set $S_{v,i} = \{v, v_0,\dots,v_{i-1},v_{i+1},\dots,v_n\}$ has $\Vol(S_{v,i}) \leq n!$ (see proof of Theorem 14.2.1 in \cite{DL}).
Hence the coefficients of $v$ are bounded by a bound depending only on $S$ and this proves that $\Hyp(S,D')$ is polyhedral. \qed

For a Delaunay simplex set $S$ of volume $1$, the cone $\Hyp(S,\ZZ^n)$ is called the
{\em hypermetric cone} and is studied in \cite{DL}.
For other simplices, they are called {\em Baranovski cone}
in \cite{bookschurmann}.
The facets of the Baranovski cones are determined up to dimension $6$ in
\cite{BaRy}.
There is a correspondence between facets of $\Hyp(S,D)$ and repartitioning sets
$P$ with $S\subset P\subset D$.
That is the inequality $f(x)\geq 0$ defines a facet of
$\Hyp(S,D)$ if and only if $S\cup\{x\}$ is a repartitioning set.

\section{Connectivity results}\label{SectionConn}

For a given Delaunay polyhedron $D$, let us write
\begin{equation*}
{\mathcal C}_{r,d}(D)=\left\{\begin{array}{l}
D'\subset D\mbox{~s.t.~} D' \mbox{~Delaunay~polyhedron}\\
\rankperf(D')=r\mbox{~and~}\dim L(D')\leq d
\end{array}\right\}.
\end{equation*}

If ${\mathcal A}$ and ${\mathcal B}$, are sets of Delaunay polyhedra then the graph $Gr({\mathcal A},{\mathcal B})$ is the graph on ${\mathcal A}$ with two Delaunay polyhedra $D_1, D_2\in {\mathcal A}$ adjacent if and only if $D_1\cap D_2\in {\mathcal B}$.

\begin{theorem}
If $D$ is a Delaunay polyhedron of perfection rank $r$ and degeneracy degree $d$ then ${\mathcal C}_{r+1, d}(D)$ is decomposed into a finite number of orbits under $\Aut(D)$.
\end{theorem}
\proof Without loss of generality, we can write $D = P + \ZZ^d$.
Let us write $D'\subset D$ as $D' = P' + L'$ with $k=\dim\,L'$.
By applying an element of $\Aut(D)$, we can assume that $L'=\ZZ^k$.
So, without loss of generality, we can assume that $L'=0$.
Let us take a Delaunay simplex set $S$ in $D'$; its volume is bounded by $\max\,PossVol(n)$.
Again by using $\Aut(D)$ we can find
a constant $C'$ such that the absolute value of the coordinates of $S$
in $P + \ZZ^d$ are bounded by $C'$. 
The polyhedrality of the cones $\Hyp(S, D)$ implies the finiteness. \qed

\begin{lemma}\label{TheCriticalLemma}
If ${\mathcal C}$ is a polyhedral cone, $F$ a face of ${\mathcal C}$ and $e$, $e'$ are two extreme rays which are not contained in $F$, then $e$ and $e'$ are connected by a path which does not intersect $F$.
\end{lemma}
\proof By taking the intersection ${\mathcal C}\cap H$ with $H$ a suitable hyperplane, we can transform ${\mathcal C}$ into a polytope $P$ and $e$, $e'$ into vertices of $P$.
We can find an affine function $\phi$ such that $\phi(x)\geq 0$ is a valid inequality on $P$ and $\phi(x)=0$ defines the face $F\cap H$ of $P$.
By maximizing the function $\phi$ over $P$ and using the simplex algorithm (see \cite{schrijver,zieglerLectureOnPolytopes}), we can find paths $p(v,v_{opt})$, $p(v',v_{opt})$ from $v$, $v'$ to an optimal vertex $v_{opt}$ such that $\phi$ is monotone on both paths.
Since $\phi(v)>0$ and $\phi(v')>0$, such paths avoid the face $F$ and put together gives the required path. \qed

\begin{theorem}\label{ErdahlConnectivityResult}
If $D$ is a Delaunay polyhedron of perfection rank $r$ and degeneracy degree $d\geq 1$ then
\begin{equation*}
Gr\left({\mathcal C}_{r+1,d}(D), {\mathcal C}_{r+2,d-1}(D)\right)
\end{equation*}
is connected.
\end{theorem}
\proof Let us take two Delaunay polyhedra $D_a$ and $D'_a$ in $D$ of perfection
rank $r+1$.
Let us take two Delaunay simplex sets $S$, $S'$ contained in $D_a$, $D'_a$.
By using Theorem \ref{ConnectivitySimplices} we can find a chain of
simplices $(S_i)_{0\leq i\leq m}$ with $S_i\subset D$, $S_0=S$ and $S_m=S'$.
Denote by $R_i=S_i\cup S_{i+1}$ the repartitioning set.
Write $e_{-1}$, $e_{m}$ to be the extreme rays in $\Hyp(S_0,D)$,
$\Hyp(S_m,D)$ corresponding to $D_a$, $D'_a$.
For each Delaunay simplex set $S_i$ we consider the cone $\Hyp(S_i,D)$.
The extreme rays correspond to Delaunay polyhedra of rank $r+1$.
If a Delaunay polyhedron $D'\subset D$ has degeneracy degree $d$ then 
necessarily $L(D')=L(D)$.
We define the restricted trace function to be
\begin{equation*}
\phi(f) = Tr\left(\left. \Quad(f)\right\vert_{L(D)}\right).
\end{equation*}
A function $f\in Erdahl_{supp}(D)$ has $Z(f)\cap D$ of degeneracy degree $d$ if
and only if $\phi(f)=0$. The hyperplane $\phi(f)=0$ determine
a face $F_i$ of the cone $\Hyp(S_i,D)$.
The intersection is
\begin{equation*}
\Hyp(S_i,D) \cap \Hyp(S_{i+1},D) = \Hyp(R_i, D).
\end{equation*}
Thus we can find a ray $e_i$ in $\Hyp(R_i,D)$, which is not contained in $F_i$.
Since $\Hyp(S_i,D)$ is polyhedral, by Lemma \ref{TheCriticalLemma} there
exists a path from $e_{i-1}$ to $e_i$ in $\Hyp(S_i,D)$ that avoids the face
$F_i$.
So, by putting all the paths together, we got the required connectivity
result. \qed

\begin{lemma}\label{DiamondOnErdahl}
If $D_1$ and $D_3$ are two Delaunay polyhedra of perfection rank $r$ and $r+2$ with $D_3\subset D_1$ then there exist exactly two Delaunay polyhedra $D_{2,1}$ and $D_{2,2}$ with
\begin{equation*}
D_3 \subset D_{2,i}\subset D_1
\end{equation*}
with $\rankperf(D_{2,i})=r+1$.
\end{lemma}
\proof Since $D_3$ is a Delaunay polyhedron there exist a Delaunay simplex set $S\subset D_3$.
The Delaunay polyhedra $D_3$, $D_1$ correspond to faces $F_3$, $F_1$ of dimension $r+2$, $r$
in the cone $\Hyp(S, \ZZ^n)$.
It is well known from polytope theory  \cite[Theorem 2.7.(iii)]{zieglerLectureOnPolytopes}
that there are exactly two faces $F_{2,1}$, $F_{2,2}$ containing $F_1$ and contained in $F_3$.
Those gives the corresponding Delaunay polyhedron. \qed

By using this theorem, we are able to compute inductively the Delaunay
polyhedra in $\ZZ^n$. The property with the degeneracy degree ensures that
we are able to effectively reduce the complexity of the computation at each
step and thus we are reduced in the end to computation with Delaunay
polyhedra of degeneracy $0$, i.e. polytopes for which polytopal methods
exist.

\section{Algorithms}\label{SectionAlgo}

In \cite{symsurvey} a general survey of methods for computing dual description
of highly symmetric polytopes with many facets are presented.
Among the method presented there, we want to adapt the Recursive Adjacency
Decomposition Method to our situation, i.e. to a case with an infinite group
and an infinity of defining inequalities.

\subsection{Computing $\Aut(D)$}\label{AlgoComputAutD}
In this subsection we explain the techniques needed to 
compute $\Aut(D)$, compute $\Stab_{D'}(\Aut(D))$ and split orbits.
In the decomposition of Theorem \ref{AutomAut_D} the only component
that is not clear is $\Aut(D_1)$, i.e. the computation of the
automorphism group of a Delaunay polytope.
For that purpose the methods of \cite{ComplexityVoronoiDSV} can be used. 
The one that we are using is the method of isometry groups.

\begin{definition}
Let $D$ be a Delaunay polyhedron. Take a lattice $L'$ with $L'\oplus_{\ZZ} L(D) = \ZZ^n$.
Take a basis $w_1, \dots, w_r$ of $L'$. Denote by $v_1, \dots, v_m$ the expression of
the vertices of $P_{L'}(D)$ in the basis $(w_i)$.\\
We define the matrix $Q$ by
\begin{equation*}
Q = \sum_{i=1}^m \left(\begin{array}{c}
1\\
v_i
\end{array}\right) (1, v_i^t).
\end{equation*}
From then we define the distance function $f_D: \ZZ^n\times \ZZ^n \mapsto \RR$ by
\begin{itemize}
\item $f_D(x, x') = \phi(x) Q^{-1} \phi(x')^T$
\item with $\phi(x) = (1, u_1, \dots, u_r)$ if 
\begin{equation*}
x = u_1 w_1 + \dots + u_r w_r + z \mbox{~and~} z\in L(D).
\end{equation*}
\end{itemize}

\end{definition}

The construction of the matrix $Q$ and its inverse above is relatively standard.
We used it first in \cite{PerfectDim8} and further work on this are done
in \cite{symsurvey,GroupPolytopeLMS}.

The interest of this construction is that it allows to compute automorphism groups.
\begin{theorem}\label{ScalProdPreservation}
  Let $D$ be a Delaunay polyhedron. It holds:

  (i) If $u\in \Aut(D)$ then we have $f_D(u(x), u(y) ) = f_D(x,y)$ for $x,y\in D$.

  (ii) If $L'$ is a sublattice such that $L'\oplus_{\ZZ} L(D) = \ZZ^n$ and $u$ is a permutation of $P_{L'}(D)$ such that $f_D(u(x), u(y) ) = f_D(x,y)$ for $x,y\in D$
  then $u$ is induced by an affine rational transformation of $L'\otimes \RR$.
\end{theorem}
\proof (i) By its construction if $z\in L(D)$ then we have $f_D(x + z,y) = f_D(x,y)$.
Thus if $u\in \Aff(D)$ then $u$ preserve $f_D$.
On the other hand if $u\in \Aut(P_{L'}(D))$ then we can see by summation
that $u$ preserve $f_D$. The proof is available for example in \cite{PerfectDim8,symsurvey}.

(ii) The reverse implication is also available from \cite{PerfectDim8,symsurvey}. \qed

Let us denote by $\Aut_{\QQ}(P(D))$ the group of rational transformations preserving
$P(D)$. By Theorem \ref{ScalProdPreservation}.(ii) we have $\Aut(P(D)) = \AGL_{r}(\ZZ) \cap \Aut_{\QQ}(P(D))$
with $r$ the dimension of $L'$.
The computation of $\Aut_{\QQ}(P(D))$ is done efficiently by using known partition backtracking
software such as \cite{nauty}; see \cite{symsurvey,GroupPolytopeLMS} for more details.
In the cases considered in this paper the number of vertices is quite small and this
computation is very easy.

A Delaunay polytope is called {\em generating} if difference between its vertices
generate $\ZZ^n$. If a Delaunay polytope is non-generating, then it is
actually a Delaunay polytopes for more than one lattice. If $P(D)$ is generating then we
have $\Aut(P(D)) = \Aut_{\QQ}(P(D))$ and we are done. Otherwise, we can apply some of the
strategies listed in \cite[Section 3.1]{GroupPolytopeLMS}. Here, the situation is particularly
simple and the simplest strategy of iterating over the group elements and keeping the integral
ones works very well.
Also note that the above methods with only slight modifications work for
testing equivalence of Delaunay polyhedra.

\subsection{Computing stabilizers}
We now give methods for computing stabilizers of Delaunay polyhedra, more precisely
of the transformations preserving two polyhedra $D\subset D'$ which occurs in our computations.

Let us select a lattice $L'$ such that $L'\oplus_{\ZZ} L(D) = \ZZ^n$.
Denote by $G_1$ the group occurring in Theorem \ref{InclusionAff}.

Let us define the following function on $P_{L'}(D)$ by
\begin{equation*}
f_{D,D'}(x, y) = \left( f_D(x,y), f_{D'}(x,y) \right).
\end{equation*}
By Theorem \ref{ScalProdPreservation} the elements of the group $G_1$
must preserve the function $f_{D,D'}$. Thus we may use the partition backtrack algorithm
of Section \ref{AlgoComputAutD} to get the group $\Aut(f_{D,D'})$
of permutations preserving $f_{D,D'}$.

Then we obtain the group $G_1$ by keeping only the elements that are in $\AGL_n(\ZZ)$ and preserve $D'$.
This is possible since the group $\Aut(f_{D,D'})$ is finite and of moderate size in most cases.
All the algorithms above have equivalents for testing equivalence and of course what has been
done for pairs $D\subset D'$ of Delaunay polyhedra can be extended to triples $D\subset D'\subset D''$.

\subsection{Splitting orbits}
Suppose that we have an orbit $Gx$ of an element $x$ under a group $G$.
For a subgroup $H\subset G$ we wish to decompose $G x$ into orbits $H x_i$.
Such an orbit splitting decomposition
\begin{equation*}
G x = \cup_{i=1}^m H x_i \mbox{~with~}x_i=g_i x
\end{equation*}
is equivalent to a double coset decomposition 
\begin{equation*}
G = \cup_{i=1}^m H g_i \Stab_G(x).
\end{equation*}

In the case of interest to us we have $G=\Aut(D)$, $H=\Aut(D)\cap \Aut(D')$ for $D$, 
$D'$ Delaunay polyhedra with $D\subset D'$ and $x$ a Delaunay polyhedron included in $D$.
Since a priori $\Aut(D)$ is infinite we cannot apply standard tools from computer algebra
software such as {\tt GAP} \cite{gap}.
By the finiteness result Theorem \ref{UniqueDecomposition}.(ii) we can find a coset decomposition
\begin{equation*}
G = \cup_{i=1}^m H g_i.
\end{equation*}
However, it is not a double coset decomposition, i.e. we can have 
$H g_i \not= H g_j$ but still have $Hg_i x = Hg_j x$.
Therefore, we need to eliminate duplicate in order to do the orbit splitting.

\subsection{The flipping algorithm}
Suppose that $D_1\subset D_2\subset D_3$ are Delaunay polyhedra having
\begin{equation*}
\rankperf\, D_1=1+\rankperf\, D_2=2+\rankperf\, D_3.
\end{equation*}
By Lemma \ref{DiamondOnErdahl} we know that there exist a unique 
Delaunay polyhedron $D'_2$ having $D_1\subset D'_2\subset D_3$
and $D_2\not= D'_2$.

We can find functions $f_i\in Erdahl(n)$ such that $Z(f_i)=D_i$.
We can also assume that $C_{D_i}= \RR f_i\oplus C_{D_{i+1}}$ for $i=1$, $2$.
We need to find $f'_2\in Erdahl(n)$ such that $Z(f'_2)=D'_2$.

If $L(D_3)=0$ then $D_3$ is a polytope, i.e. it has a finite number
of vertices and the algorithm is called the {\em gift wrapping} procedure
\cite{symsurvey,PerfectDim8,CR}.
If $L(D_3)\not= 0$ we have to modify the algorithm in order to take care of the
fact that we have an infinity of vertices by writing an iterative
algorithm. This is quite similar to the flipping in the Voronoi
algorithm \cite{bookschurmann}.

\begin{algorithm}\label{Flipping_algorithm}
\KwData{Delaunay polyhedra $D_1$, $D_2$, $D_3$ with $D_i=Z(f_i)$, $f_i\in Erdahl(n)$, $D_1\subset D_2\subset D_3$ and
\begin{equation*}
\rankperf\,D_1-2=\rankperf\, D_2-1=\rankperf\,D_3.
\end{equation*}}
\KwResult{$f'_2\in Erdahl(n)$, Delaunay polyhedra $D'_2=Z(f'_2)$,  $D'_2\not= D_2$ and $D_1\subset D'_2\subset D_3$}
${\mathcal V} \leftarrow \emptyset$\\
\Repeat{${\mathcal L}$ has rank $2$}{
  $v \leftarrow \mbox{random~element~of~}Z(f_3)$\\
  ${\mathcal V} \leftarrow {\mathcal V}\cup \{v\}$\\
  ${\mathcal L} \leftarrow (f_1(w), f_2(w))\mbox{~for~}w\in {\mathcal V}$
}
\Repeat{$f'_2\geq 0$ on $Z(f_3)$}{
  $\{(\alpha_1, \beta_1), (\alpha_2,\beta_2)\} \leftarrow \mbox{~generators~of~extreme~rays~of~}{\mathcal L}$.\\
  $f'_2\leftarrow \{\alpha_1 f_1 + \beta_1 f_2, \alpha_2 f_1 + \beta_2 f_2\} - \RR f_2$.\\
  \If{there is a $v\in Z(f_3)$ with $f'_2(v) < 0$}{
    ${\mathcal V} \leftarrow {\mathcal V}\cup \{v\}$.\\
    ${\mathcal L} \leftarrow (f_1(w), f_2(w))\mbox{~for~}w\in {\mathcal V}$
  }
}
\Repeat{$f'_2\geq 0$ on $\ZZ^n$}{
  $f'_2\leftarrow f'_2 + f_3$\\
}
$f'_2\leftarrow f'_2 + f_3$
\caption{Flipping algorithm}
\end{algorithm}
The non-negativity test for $f'_2$ on $\ZZ^n$ is done by solving a closest vector problem.
The non-negativity test on $D_3$ is done by decomposing it into $\{v_1 + L(D_3)\} \cup \dots \cup \{v_m + L(D_3) \}$.
The non-negativity is tested by $m$ closest vector problems.
The final operation on $f'_2$ is done to ensure that $f'_2(x)>0$ if $x\notin Z(f_3)$.

\subsection{The recursive adjacency decomposition method}

Given a Delaunay polyhedron $D$ of perfection rank $r$ the algorithm of this section will give the
orbits of Delaunay polyhedra $D'\subset D$ of perfection rank $r+1$.
If $\degrk(D)=0$ then the computation of the orbits of Delaunay polyhedra can be achieved by
Theorem \ref{BasicResultErdahlCone}.(iii).

If $\degrk(D)>0$ we have to proceed differently. By Theorem \ref{ErdahlConnectivityResult}
we can limit ourselves to Delaunay polyhedra with $\degrk(D)\leq \degrk(D)-1$.
The algorithm takes one initial Delaunay polyhedron of perfection rank $r+1$ and computes
the adjacent Delaunay polyhedron of perfection rank $r+1$.
If an obtained Delaunay polyhedron is not equivalent to an existing one then we insert it into the
list.
We iterate until all orbits have been treated.
The computation of the adjacent Delaunay polyhedra adjacent to a Delaunay polyhedron $D'$
requires the computation of orbits the Delaunay polyhedra contained in $D'$.
Thus we have a recursive call to the algorithm. Fortunately the degeneracy degree diminish
by at least $1$ so there is no infinite recursion.

The mapping from $\MF_1$ to $\MF_2$ is done using the orbit splitting
procedure.
Finding $D''$ from $D_2$, $D'$ and $D$ is done using the flipping procedure.

There is a degree of choice in the initial Delaunay polyhedron $D_{init}$.
The standard choice is if $D=P(D) + L(D)$ with $L(D)=\ZZ v_1 + \dots + \ZZ v_{\degrk(D)}$ to take $D_{init}$ a Delaunay polyhedron of the form
\begin{equation*}
D_{init} = P(D) + \{0, v_1\} + \ZZ v_2 + \dots + \ZZ v_{\degrk(D)} .
\end{equation*}
The schematic of the algorithm is shown in Algorithm \ref{Enumeration_inequivalent}.

\begin{algorithm}\label{Enumeration_inequivalent}
\KwData{Delaunay polyhedron $D$ of perfection rank $r$}
\KwResult{Set $\MF=\{D_1, \dots, D_m\}$ of Delaunay polyhedron of perfection rank $r+1$ inequivalent under $\Aut\, D$}
\eIf{$\degrk(D)=0$}{
  $U \leftarrow$ orbits of facets of $Erdahl_{supp}^*(D)$.\\
  $\MF \leftarrow$ orbits of Delaunay polyhedron from Theorem \ref{BasicResultErdahlCone}.(iii).
}{
  $T \leftarrow \{D_{init}\}$ with $D_{init}\in Erdahl(n)$, $\rankperf\,Z(D_{init})=r+1$ and $\degrk(D_{init})=\degrk(D)-1$.\\
  $\MF \leftarrow \emptyset$.\\
  \While{there is a $D' \in T$ with $\degrk(D')\leq \degrk(D)-1$}{
    $\MF \leftarrow \MF \cup \{D'\}$.\\
    $T \leftarrow T \setminus \{D'\}$.\\
    $\MF_1 \leftarrow$ orbits of Delaunay polyhedra of perfection rank $r+2$ in $D'$ under $\Aut\, D'$.\\
    $\MF_2 \leftarrow$ orbits of Delaunay polyhedra of perfection rank $r+2$ in $D'$ under $\Stab(D', D)$.\\
    \For{$D_2 \in \MF_2$}{
      find Delaunay polyhedron $D''\subset D$ of perfection rank $r+1$ with $D_2 = D' \cap D''$.\\
      \If{$D''$ is not equivalent under $\Aut\, D$ to an element of $\MF\cup T$}{
        \eIf{$\degrk(D'') = \degrk(D)$}{
          $\MF \leftarrow \MF \cup \{f'\}$
        }{
          $T \leftarrow T \cup \{f'\}$
        }
      }
    }
  }
}
\caption{Enumeration of inequivalent sub Delaunay polyhedra}
\end{algorithm}

\section{Perfect Delaunay polytopes in dimension $7$}\label{DimSeven_perfectDelaunay}

In the enumeration of inhomogeneous perfect form in dimension $7$ we need to
describe the Delaunay polytopes that will occur.
The list of perfect Delaunay polyhedra in dimension $7$ is thus
\begin{enumerate}
\item $\{0,1\}\times \ZZ^6$.
\item $2_{21} \times \ZZ$ with $2_{21}$ the Schl\"afli polytope
\item $3_{21}$ the Gosset polytope \cite{DL}.
\item $ER_7$ the polytopes discovered by Erdahl and Rybnikov \cite{35tope,InfiniteER2002}.
\end{enumerate}
The geometry of the Schl\"afli and Gosset polytope are described in more details in \cite{DL,PerfectDelaunayLowDim}.

An {\em affine basis} of a $n$-dimensional Delaunay polytope $D$
is a family of $n+1$ vertices $v_0$, \dots, $v_n$ such that for any
vertex $v$ of $D$ there exist $\lambda_i\in \ZZ$ such that
$v=\sum_{i=0}^n \lambda_i v_i$. The perfect Delaunay polytopes of
dimension $7$ have an affine basis but it is possible that in higher
dimension there are perfect Delaunay polytopes without affine basis.
It is known that in dimension at least $12$, there are Delaunay polytopes
with no affine basis \cite{RankComput}.
Also the perfect Delaunay polytopes of dimension $7$ are generating.
Note that in \cite{SecInfinite} we found some non-generating perfect Delaunay polytopes
for $n\geq 13$.
We have $\rankperf(\{0,1\}^n)=n$ (see for example \cite{DelaunaySix,DL}).

In terms of computation, the overwhelming majority of the time is spent computing
the rank $2$ faces of $\{0,1\}\times \ZZ^6$.
By the recursive approach chosen, the method requires the computation of the facets
of $Erdahl_{supp}^*(\{0,1\}^7)$ and so by Theorem \ref{CUTn_and_ErdahlHn} of the
facets of $\CUTP_8$.
We actually computed the list of orbits of facets of $\CUTP_8$ (and some other graph
cut polytopes) in \cite{CUT8_facet}.
In dimension $8$ the partial enumeration algorithm of
\cite{ANewAlgorithm} found $27$ perfect Delaunay polytopes and it is likely
that the list is complete. But to prove its completeness by using the
method of this work would require the determination of all facets of $\CUTP_9$
and this is very hard \cite{CR}.
In \cite{ANewAlgorithm} a partial enumeration of perfect Delaunay polytopes
was done with only Delaunay polyhedra with $L(D)=0$ being considered.
The two perfect Delaunay polytopes of dimension $7$ were determined in this work
and our enumeration proves that the list is complete.

The implementation is available from \cite{Polyhedral} and uses the {\tt GAP}
computer algebra system \cite{gap}.

\section{Classification of Delaunay simplices in dimension $7$}\label{DimSeven_simplices}

Formula \eqref{UpperBound_PossVol} gives $187$ as an upper bound on the volume
of Delaunay simplex sets.
With this upper bound we can devise an algorithm for enumeration of Delaunay simplex sets,
which will unfortunately prove inefficient:

\begin{remark}\label{ExpansionAlgorithm}
Suppose we have a list of types of Delaunay simplex sets in dimension $n-1$.
If $S = \{v_0, \dots, v_{n-1}, v_n\}$ is a Delaunay simplex set of dimension $n$
then $\{v_0, \dots, v_{n-1}\}$ is a $n-1$ dimensional Delaunay simplex sets of the
lattice $v_1 - v_0$, \dots, $v_{n-1} - v_0$.

If we have a $n-1$ dimensional Delaunay simplex set $v_0=0$, $v_i\in \ZZ^{n-1}$ for
$1\leq i\leq n-1$ and is of volume $v$.
We write the $n$-dimensional simplex as $v'_0 = (v_0, 0)$, $v'_1 = (v_1, 0)$, \dots, $v'_{n-1} = (v_{n-1}, 0)$
and $v'_n = (x_1, \dots, x_{n-1}, x_n)$ with $x_n > 0$.
The volume of the simplex defined by $(v'_i)$ is $v x_n$.

For a fixed $x_n$ the number of possibilities for $(x_1, \dots, x_{n-1})$ is $x_n^{n-1} / v$.
\end{remark}

The Delaunay simplices of dimension $6$ were classified in \cite{Ba} and so we could apply
the algorithm of Remark \ref{ExpansionAlgorithm}.
Unfortunately the number of possibilities to be applied is very large, of the order
of $187^6$, on which we have to apply Algorithm \ref{TestRealizability}.

Therefore, we need a different method:

\begin{lemma}\label{Simplices_NotInPerfect}
Let $S$ be a Delaunay simplex set that is not contained in any perfect Delaunay
polyhedron different from $\{0,1\}\times \ZZ^{n-1}$. The possibilities are:
\begin{enumerate}
\item For $n\leq 4$ the Delaunay simplices of volume $1$.
\item For $n=5$ the Delaunay simplex of volume $1$ or $2$.
\item For $n=6$ or $7$ there are no possibilities.
\end{enumerate}

\end{lemma}
\proof Let us take a Delaunay simplex set $S = \{v_0, \dots, v_n\}$.
We can assume that $v_0$ is located at the origin by using translation if necessary.
For each $1\leq i\leq n$ let us define $\ell_i$ the linear form on $\RR^n$ such that $\ell_i(v_i)=1$ and $\ell_i(v_j)=0$ for $i\not=i$.
Any Delaunay polyhedron $D$ isomorphic to $\{0,1\}\times \ZZ^{n-1}$ and such that $S \subset D$ corresponds to a linear form $\ell$ on $\RR^n$ such that $\ell(v_i)\in \{0,1\}$ and
\begin{equation*}
D=\{x\in \ZZ^n \mbox{~$\vert$~} \ell(x) = 0\mbox{~or~} 1\}.
\end{equation*}
The linear form $\ell$ is then called {\em admissible} and the corresponding quadratic function is $q_{\ell}(x)= \ell(x)( \ell(x) - 1)$.

Let us denote by $S\subset \{1,\dots,n\}$ the set of points $i$ such that $\ell(x)=1$.
Clearly, one can write
\begin{equation*}
\ell = \sum_{i\in S} \ell_i.
\end{equation*}
A function $\ell$ is admissible if and only if $\ell$ is integral valued on $\ZZ^n$.
If it is not integral valued then there exists a $v\in \ZZ^n$ such that $0 < \ell(v) < 1$
which implies that $q_{\ell}(v) < 0$ which is not allowed.
If it is integral valued then $D=Z(q_{\ell})$ is equivalent to $\{0,1\}\times \ZZ^{n-1}$.
In the following, for a set $S\subset \{1, \dots, n\}$ we write
\begin{equation*}
\ell_S = \sum_{i\in S} \ell_i \mbox{~and~} v_S = \sum_{i\in S} v_i.
\end{equation*}
Let us define
\begin{equation*}
{\mathcal S} = \left\{ S\subset \{1, \dots, n\} \mbox{~s.t.~} \ell_S\mbox{~is~integral~valued}\right\} .
\end{equation*}
Let us denote by $\ZZ {\mathcal S}$ the $\ZZ$-span of the elements $v_S$ for $S\in {\mathcal S}$. This defines a lattice ${\mathcal L}$ of $\ZZ^n$. 
Since $S$ is contained only in Delaunay polyhedra isomorphic to $\{0,1\}\times \ZZ^{n-1}$ the set of the function $q_{\ell_S}$ is full-dimensional. This implies that $\left\vert{\mathcal S}\right\vert \geq n(n+1)/2$ and that the lattice ${\mathcal L}$ is actually full dimensional. Denote by $h$ its index.

The set of the function $\ell_S$ is also full-dimensional in $(\ZZ^n)^*$. Its index is
\begin{equation*}
h/\Vol(S)\geq 1.
\end{equation*}
We are interested in the point sets of the form $\{0,1\}^n\cap L$ with $L$ an affine subspace of $\ZZ^n$. 
By direct enumeration we obtain the full list of $3363$ orbits of such points for $n=7$.
By selecting the point sets whose cone of functions $q_{\ell_S}$ is full-dimensional we get
an upper bound of $3$ on the index $h$ and so an upper bound of $3$ on the possible volumes
of such simplices.

With volume at most $3$ we can apply the algorithm implied by Remark \ref{ExpansionAlgorithm}.
Each facet of such a Delaunay simplex is also a Delaunay simplex of one dimension lower.
Therefore, we can use previous enumeration result to get a list of $796$ possible candidates
of $7$-dimensional Delaunay simplices.
We then use Algorithm \ref{TestRealizability} for checking which ones of them are indeed
Delaunay simplices. This gives $6$ cases (the ones of Table \ref{TableFundamentalSimplices}
of volume at most $3$). Each one of them is also contained in a Delaunay polytope $ER_7$
and so there is no such Delaunay simplices in dimension $7$.

Dimension $n\leq 6$ follows from known results. \qed

{\bf Proof of Theorem \ref{EnumerationSimplices}}:
If $S$ is a Delaunay simplex set then $\Hyp(S,\ZZ^7)$ is a full-dimensional polyhedral cone,
i.e. defined by a finite number of inequalities and having a finite number of extreme rays.
Any such extreme ray corresponds to a perfect Delaunay polyhedron $D$.
We have $|S|=8$ and $S$ defines a face of the cone $Erdahl_{supp}^*(D)$.

By Lemma \ref{Simplices_NotInPerfect}, $S$ has to be contained in 
a Delaunay polyhedron of type $3_{21}$, $ER_7$ or $2_{21}\times \ZZ$.

The perfect Delaunay polyhedron $ER_7$ has $35$ vertices which matches
the lower bound given by Proposition \ref{LowerBoundNrVertices}.
As a consequence any $8$-element subset of $ER_7$ defines a face of
$Erdahl_{supp}^*(ER_7)$.
The automorphism group of $ER_7$ has size $1440$ and by using it one
can get easily the $9434$ orbits of $8$-element subsets of $ER_7$.
Actually all $11$ types of simplices occur this way.

The Gosset polytope $3_{21}$ has $56$ vertices and the automorphism group
is equal to the Weyl group of the root lattice $\mathsf{E}_7$. We found
$521$ orbits of $8$-element sets in $3_{21}$, $474$ of them correspond
to faces of $Erdahl_{supp}^*(3_{21})$.

For the perfect Delaunay polyhedron $2_{21}\times \ZZ$ we have to
proceed differently since the number of points to be considered is infinite.
We have to enumerate the possible $8$-point subsets of $2_{21}\times \ZZ$ of
volume at most $187$ (Formula \eqref{UpperBound_PossVol}) up to the action of $\Aut(2_{21}\times \ZZ)$.
The $8$ points are expressed in the form $v_i=(w_i, h_i)$ with $w_i\in 2_{21}$
and $h_i\in \ZZ$.
The set of points $(w_i)_{1\leq i\leq 8}$ must define a $6$-dimensional affine space.
Thus $7$ of them, say $(w_i)_{1\leq i\leq 7}$, must be sufficient to define a $6$-dimensional
Delaunay simplex set $S_{Sch}$.
An exhaustive enumeration on the $27$ vertices of $2_{21}$ gives $31$ types up
to isomorphism.
The volume $\Vol(S_{Sch})$ can be $1$, $2$ or $3$. If the volume is $1$, then we can use an element
of $\Aff(2_{21}\times \ZZ)$ and obtain $h_i=0$ for $1\leq i\leq 7$.
For higher volumes, the situation is more complicate but by using $\Aff(2_{21}\times \ZZ)$
and linear algebra we we can reduce to $\Vol(S_{Sch})$ possibilities, i.e. $2$ or $3$.
For the last point $(v_8, h_8)$ we have $27$ possibilities for $v_8$ and a finite number for
$h_8$ due to the upper bound of $187$.
We then apply Algorithm \ref{TestRealizability} to test realizability of the finite list
of possible cases. This gives us the $11$ possible simplices. \qed

\section*{Acknowledgments}
The author thanks Viacheslav Grishukhin and Achill Sch\"urmann for useful discussion on this work.

\bibliographystyle{amsplain_initials_eprint}
\bibliography{LatticeRef}

\end{document}